\documentclass{ifacconf}

\usepackage{algorithmic}
\usepackage{amsmath,amssymb}
\usepackage{color}
\usepackage{bm}
\usepackage{subfigure}
\usepackage{graphicx}      
\usepackage{natbib}        
\usepackage{booktabs}
\usepackage{multirow}
\usepackage{upgreek}
\usepackage{overpic}
\usepackage{enumitem}

\usepackage{xspace} 
\usepackage[
textsize=scriptsize
]{todonotes}

\newtheorem{theorem}{\textbf{Theorem}}

\definecolor{dblue}{rgb}{ 0.00, 0.00, 0.60 }

\begin{document}

\begin{frontmatter}

\title{Interval Superposition Arithmetic for Guaranteed Parameter Estimation}

\author[ST]{Junyan Su}
\author[ST]{Yanlin Zha}
\author[ST]{Kai Wang}
\author[ST]{Mario E. Villanueva}
\author[STUBA]{Radoslav Paulen}
\author[ST]{Boris Houska}

\address[ST]{School of Information Science and Technology, ShanghaiTech
University, 319 Yueyang Road, Shanghai 200031, China. ({\tt \{zhayl,wangk,meduardov, borish\}@shanghaitech.edu.cn)}}
\address[STUBA]{Faculty of Chemical and Food Technology, Slovak University
of Technology in Bratislava, Radlinskeho 9, Bratislava, Slovakia.
({\tt radoslav.paulen@stuba.sk})}

\begin{abstract}                
The problem of guaranteed parameter estimation (GPE) consists in enclosing the
set of all possible parameter values, such that the model predictions match the
corresponding measurements within prescribed error bounds. One of the bottlenecks
in GPE algorithms is the construction of enclosures for the image-set of factorable
functions. In this paper, we introduce a novel set-based computing method called
interval superposition arithmetics (ISA) for the construction of enclosures of
such image sets and its use in GPE algorithms. The main benefits of using
ISA in the context of GPE lie in the improvement of enclosure accuracy and in the
implied reduction of number set-membership tests of the set-inversion algorithm.
\end{abstract}


\end{frontmatter}


\section{Introduction}
\label{sec::intro}

In science and engineering, the behavior of processes and systems is often
described using a mathematical model. Mathematical model development often follows
three steps: model structure specification, design (and realization) of experiments,
and estimation of unknown model parameters~\citep{fra08}. In the last step,
parameters are sought for which the model outputs match available
measurements~\citep{lju99}.

One possible way of addressing the parameter estimation problem is the use of
set-membership estimation~\citep{mil91,bai95}, also called guaranteed parameter
estimation (GPE). The GPE problem can be formulated as an identification of the
set of all possible model parameter values which are not falsified by the plant
measurements, within some prescribed error bounds. A set-inversion
algorithm~\citep[e.g. SIVIA by][]{jau93} can be applied to find such set for
nonlinear models. Here, the parameter set is successively partitioned into
smaller boxes and using exclusion tests some of these boxes are eliminated, until
a desired approximation is achieved. Since its advent, GPE has found various
applications~\citep[see e.g.,][]{mar00,jau02,lin07,has15,pau15-ima}.

The complexity of the search procedure in SIVIA is proportional to the tightness
of the interval enclosures. Considerable effort has then been invested into
developing different set-arithmetics~\citep[][such as Taylor models]{makino1996,
pau15-ima} to produce tighter enclosures of image-set of nonlinear factorable
functions. These techniques usually require computing and storing quantities such
as sensitivity information.

Here, we propose an attempt to improve GPE algorithms using a novel non-convex
set-arithmetic called Interval Superposition Arithmetic (ISA). This arithmetic
operates over Interval Superposition models (ISM), representing a piecewise
constant enclosure over a grid of the domain. Unlike a naive application of
interval arithmetic (IA) over the grid, the computational and storage complexity of
ISA is polynomial. Furthermore, it is able to exploit separable structures in the
computational graph of a factorable function. Finally, unlike Taylor model
arithmetics---which are based on local information---ISA is based on globally
valid algebraic relations. As a result, ISMs are tighter than Taylor models---at
least over large domains.

The rest of the paper is organized as follows, Section~\ref{sec::GPE}
reviews GPE and set-inversion. Section~\ref{sec::ISA} presents an
overview of ISA. An algorithm for intersecting ISMs with an interval---which
forms the basis for a set-inversion algorithm---is presented in
Section~\ref{sec::SIVISA}. It is important to notice that
the intersection algorithm runs in polynomial time, but the complexity
of computing an arbitrarily close approximation of the parameter set
is exponential. The application of the proposed algorithm to a simple
case study is shown in Section~\ref{sec::case}. Section~\ref{sec::conclusion}
concludes the paper.

\vspace{-1em}


\paragraph*{Notation} The set of real valued compact interval
vectors is denoted by $\mathbb I^{n} = \{[a,b]\subset \mathbb R^{n} \; | \;
a,b\in\mathbb R^{n},\, a\leq b\}$. Let $I=[a,b]\in\mathbb I$ and $c\in\mathbb R$,
$c+I=I+c$ we have $[a+c,b+c]$. Similarly, $cI=Ic$ denotes $[ca,cb]$ if
$c\geq 0$ ($[cb,ca]$ if $c<0$). The diameter of $I$ is denoted by
$\operatorname{diam}(I)=b-a$. Interval operations are evaluated by
IA~\citep{Moore2009}, e.g.,
\begin{align*}
[a,b] + [c,d] &= [a+b,c+d]\;,\\
[a,b] * [c,d] &= [\min\{ac,ad,bc,bd\},\max\{ac,ad,bc,bd\}]\\
\exp([a,b]) & = [\exp(a),\exp(b)]
\end{align*}


\section{Guaranteed parameter estimation}
\label{sec::GPE}

We consider a system represented by the algebraic model
\begin{equation}
\label{eq::model}
y = f(x)\;.
\end{equation}
Here, $x\in\mathbb{R}^{n_x}$ denotes unknown parameter while $y\in\mathbb{R}^{n_y}$
the (observed) output variables. The model is described by the, posibly nonlinear,
function $f:\mathbb{R}^{n_x}\to\mathbb{R}^{n_y}$.

Given $n_m\in\mathbb N$ measurements, $y^{\rm m}_{1},\ldots,y^{\rm m}_{n_m}\in
\mathbb{R}^{n_y}$, the GPE paradigm works under the assumption that true system
outputs $y^{\circ}_{1},\ldots,y^{\circ}_{N}$ can be observed only within some
bounded measurement bounds. Thus, for each $i\in\{1,\ldots,n_m\}$, we have
\begin{equation}
y^{\circ}_{i}\in y^{\rm m}_{i}+[-\eta_{i},\eta_{i}]=:Y_{i}\in\mathbb I^{n_y}\,
\end{equation}
with $\eta_{1},\ldots,\eta_{n_m}\geq 0$. The aim of GPE is to compute the set
\begin{equation}
\label{eq::gpeset}
X_{\rm e} := \left\{x\in X_{0}\ |\ \forall i\in\{1,\ldots,N\}: f(x)\in Y_{i}
\right\}\;,
\end{equation}
i.e., the set of parameters (within some admissible domain $X_{0}\in
\mathbb{I}^{n_x}$) for which the model outputs are consistent with all the
uncertain observations $Y_{i}$.

Computing~\eqref{eq::gpeset} requires intersecting the preimage of $Y_{i}$
under $f$, with the initial parameter domain, i.e.,
\begin{equation}
X_{\rm e} = \left(\bigcap_{i=1}^{n_{m}} f^{-1}(Y_{i}) \right) \cap X_{0} \;.
\end{equation}
This problem is intractable, in all but the simplest of cases, and thus
one has to settle for approximations of this set. State-of-the-art algorithms for
set inversion provide inner ($\mathbb X_{\rm int}$) and boundary
($\mathbb X_{\rm bnd}$) subpavings, i.e. lists of non overlapping interval
vectors, satisfying
\begin{equation}
\bigcup_{X\in\mathbb X_{\rm int}} X \subseteq X_{\rm e}
\subseteq \left( \bigcup_{X\in\mathbb X_{\rm int}} X \right) \cup
\left( \bigcup_{X\in\mathbb X_{\rm bnd}} X \right)\;.
\end{equation}

In a nutshell, these algorithms work by subdividing the parameter domain $X_{0}$
into smaller boxes such that $X_{0} = \bigcup_{j} X_{j}$. Set arithmetics
are then used to construct enclosures of $f$ on $X_{j}$, i.e. sets
$\overline{Y}_{j}\subset
\mathbb R^{n_y}$ satisfying
\begin{equation}
\label{eq::enclosure}
\overline{Y}_{j}\supseteq \{ f(x) \ |\ x\in X_{j} \}\;.
\end{equation}
Using the information provided by the enclosure $Y_{j}$, the following set
membership tests can be performed to classify the parameter boxes $X_{j}$ as
interior or boundary boxes:
\begin{enumerate}
\item If $\overline{Y}_{j}\subseteq Y_{i}$ for all $i\in\{1,\ldots,n_m\}$,
$X_{j}\in\mathbb{X}_{\rm int}$
\item Else, if $Y_{i}\cap f(X) = \emptyset$ for some $i\in\{1,\ldots,n_m\}$,
$X_{j}\cap X_{\rm e}= \emptyset$
\item Else, $X\in\mathbb{X}_{\rm bnd}$.
\end{enumerate}

\begin{figure}
\begin{center}
\begin{overpic}[width=0.75\columnwidth]{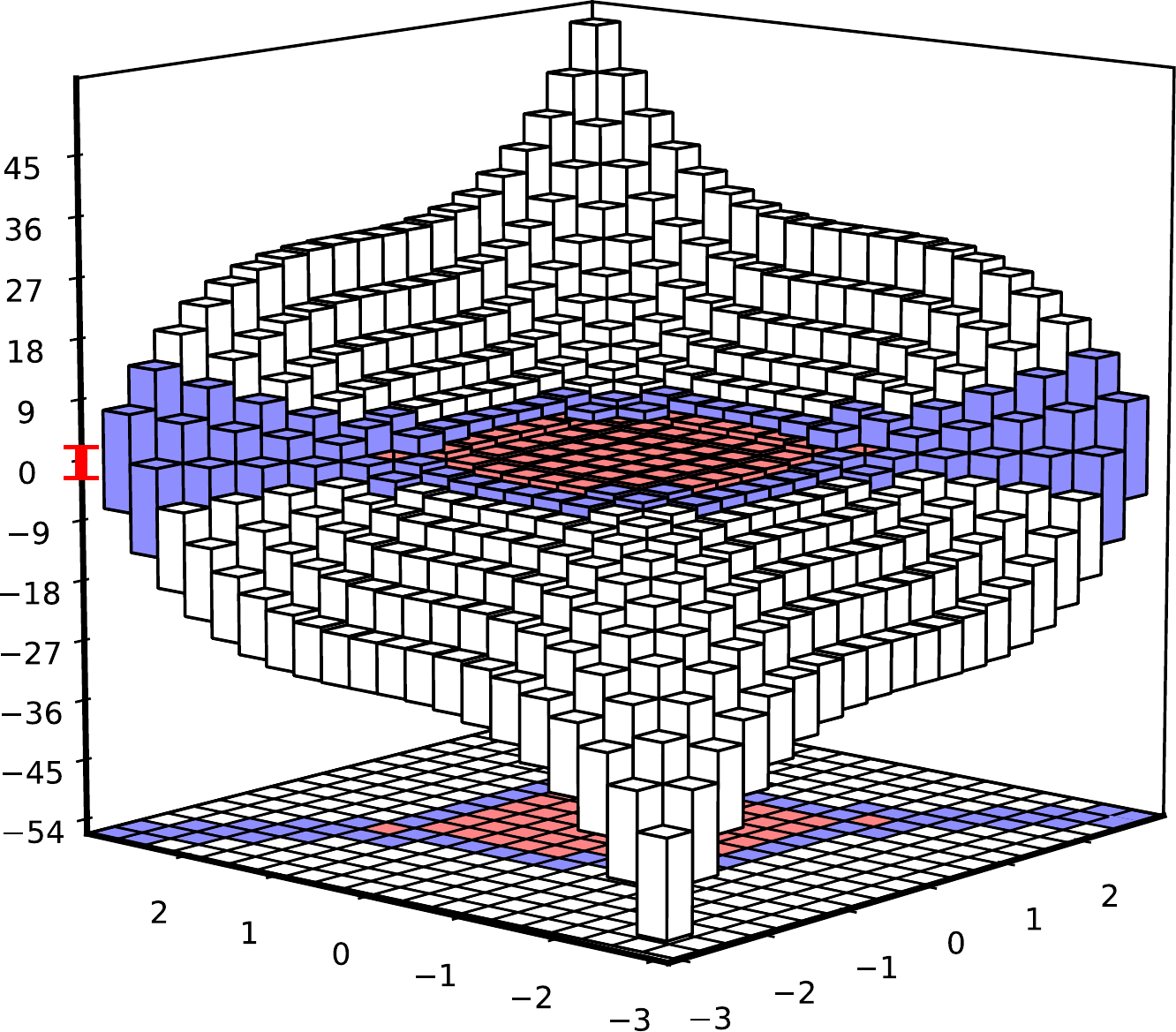}
\put(-5,49){\small $y$}
\put(80,0){\small $x_{1}$}
\put(25,0){\small $x_{2}$}
\end{overpic}
\vspace{1em}
\caption{\label{fig::sinv} Graph of an enclosure of $f=x_{1}^{3} + x_{2}^{3}$ over
$X_{0} = [-3,3]^2$ (gridded using $N=20$ subintervals at each coordinate). The
sets $\mathbb{X}_{\rm bnd}$ (blue) and $\mathbb{X}_{\rm int}$ (red) were computed
using $Y=[-2,2]$.}
\end{center}
\end{figure}

Figure~\ref{fig::sinv} shows the result of the above process for the function
$f=x_{1}^{3} + x_{2}^{3}$ over $X_{0} = [-3,3]^2$, with $Y = [-2,2]$. The
set $X_0$ has been divided into $N=20$ equidistant pieces along
each coordinate, resulting in 400 interval vectors $X_{j}$. The plot shows the
set $\bigcup_{i=1}^{N^{n_x}} \left(X_{j} \times \overline{Y}_{j}\right)$,
and its projection onto the $(x_1,x_2)$-space. The red and blue boxes
belong to $\mathbb X_{\rm int}$ and $\mathbb X_{\rm bnd}$ respectively.

In practice, the domain $X_{0}$ is subdivided iteratively by bisecting boundary
boxes, starting with $\mathbb X_{\rm bnd} = X_{0}$ and $\mathbb X_{\rm int}=\emptyset$.
The bounding, set-membership, and bisection operations are
repeated until a termination criterion, e.g.
\begin{equation}
\label{eq::maxdiam}
\forall X\in\mathbb{X}_{\rm bnd}:\quad \operatorname{diam}(X)\leq \epsilon\;,
\end{equation}
for a user-defined tolerance $\epsilon>0$, is met.
0

One of the bottlenecks of set inversion algorithms
is the over-conservatism of existing set-arithmetics, particularly over large
domains. Hence we propose to approach this problem within a novel
set-arithmetics paradigm.


\section{Interval superposition arithmetic}
\label{sec::ISA}

Interval superposition arithmetic is a novel enclosure method for the image
set of nonlinear factorable functions. It propagates nonconvex sets, called
interval superposition models, through the computational graph of the function.
Unlike Taylor~\citep{makino1996} and Chebyshev models~\citep{bat04,raj17},
ISA does not rely on local approximation methods, instead relying on global
algebraic properties and partially separable structures within the function.

\subsection{Interval superposition models}

Consider an interval domain $X=\left[\underline{x}_1,\overline{x}_{1}\right]\times
\ldots\times\left[\underline{x}_{n_x},\overline{x}_{n_x}\right]$. Now,
take a partition of $X$ into intervals of the form
\begin{equation}
\label{eq::branching}
X_{i}^{j} = [ \underline{x}_{i}+(j-1)h_i,\underline{x}_{i}+jh_{i}] \quad\text{with}
\quad h_i=\frac{\overline{x}_i-\underline{x}_i}{N}\;,
\end{equation}
for all $i\in\{1,\ldots,n_x\}$ and all $j\in\{1,\ldots,N\}$, with $N$ being a
user-specified integer.
An interval superposition model of a real-valued function $f:\mathbb{R}^{n_x}\to\mathbb{R}$
on $X$ is an interval valued function
$\Gamma:X\times\mathbb{I}^{n_x\times N}\times\mathbb{I}^{n_x}\to\mathbb{I}$, given
by
\begin{equation}
\label{eq::ISM}
\Gamma(x,A,X) = \sum^{n_x}_{i=1}\sum^{N}_{j=1} A_{i}^{j}\varphi_{i}^{j}(x)\;,
\end{equation}
with
\begin{equation}
\varphi_i^{j}(x) =
\begin{cases}
1 &\text{if} \ x_{i}\in X_{i}^{j}, \\
0 &\text{otherwise.}
\end{cases}
\end{equation}
Here, $A_{i}^{j} = \left[ \underline{A}^{i}_{j},\overline{A}^{i}_{j}\right]$ are
the components of a matrix
\begin{equation}
A =
\begin{pmatrix}
A^{1}_{1} & \ldots & A^{N}_{1} \\
\vdots & \ddots & \vdots \\
A^{1}_{n_x} & \ldots & A^{N}_{n_x}
\end{pmatrix} \in \mathbb{I}^{n_x\times N}\;,
\end{equation}
which, for a fixed $X$, completely determines the enclosure function of $f$.
Note that ISMs for functions $f:\mathbb{R}^{n_x}\to\mathbb{R}^{n_y}$
are defined by stacking ISMs for each $f_{i}$.
The matrix $A$ is constructed such that $\Gamma(\cdot,A,X)$ is a piecewise
constant enclosure function of $f$ over $X$, i.e.
\begin{equation}
\forall x\in\mathbb{X}: \quad f(x)\in \Gamma(x,A,X)\;.
\end{equation}

The name \emph{interval superposition} is motivated by the structure of the
enclosure function: At any $x\in X^{j}_{1}\times\ldots\times X^{j}_{n_x}$,
the interval $\overline{Y} = \Gamma(x,A,X)$ is given by the Minkowski sum (or
superposition) of $n_x$ interval functions
$\sum^{N}_{j=1}A^{j}_{i}\varphi^{j}_{i}(x)$.
The separable structure of ISMs allows for a storage complexity of order
${\rm\bf O}(n_xN)$, since
only $n_{x}N$ intervals need to be stored, in the matrix $A$, to represent the
$N^{n_x}$ pieces of the enclosure. In Figure~\ref{fig::sinv} the graph of an
ISM, over a partition of $X$ (with $N=20$) is shown. Although this set consists
of $400$ interval vectors (shown in red, white and blue), only 40 intervals are
stored in the matrix $A$.

This separability also allows for the global minima and maxima of
$\Gamma(\cdot,A,X)$ over $X$,
\begin{equation*}
\label{eq::rbounder}
\lambda(A) = \sum^{n_x}_{i=1}
\underbrace{\min_{j\in\{1,\ldots,N\}} \underline{A}^{j}_i}_{=:L(A_{i})}
\ \, \text{and} \ \,
\mu(A) = \sum^{n_x}_{i=1}
\underbrace{\max_{j\in\{1,\ldots,N\}} \overline{A}^{j}_i}_{=:U(A_{i})}\;,
\end{equation*}
to be computed with a complexity of order ${\rm\bf O}(n_xN)$. The interval
$\left[\lambda(A),\mu(A)\right]$ denotes the range of ISM.

\subsection{Arithmetic rules for interval superposition models}

Interval superposition arithmetics propagates ISMs through the
computational graph of a factorable function, defined by
a finite recursive composition of atom operations from a finite library
$\mathcal L = \{ \exp, \, \sin,\, +,\, *,\,\ldots\}$.

Consider the functions $g,h:X\to\mathbb{R}$, and a (possibly bivariate) atom
operation $\alpha$. Let the interval matrices $A,B\in\mathbb{I}^{n_x \times N}$
be the respective parameters for ISMs of $g$ and $h$ over $X$. In ISA, a univariate
composition rule is a map taking $A$ as an input and returning an interval matrix
$C\in\mathbb{I}^{n_x\times N}$ parameterizin an ISM such that
\begin{equation*}
\forall x\in\mathbb{X}: \quad (\alpha \circ g)(x) \in \Gamma(x,C,X)\;.
\end{equation*}
Here, $\alpha \circ g$ denotes the composition of $\alpha$ and $g$.

Bivariate composition rules in ISA are defined analogously, with the map taking
both $A$ and $B$ as inputs. Although such maps
are specific for each atom operation $\alpha$, the main steps are
outlined in Algorithms~\ref{alg::univariate} and~\ref{alg::product} for
univariate compositions and bivariate products respectively.
The addition rule in interval superposition arithmetic is simple. An interval
superposition model of $g+h$ on $X$ is parameterized by the matrix $C=A+B$, with
the sum computed componentwise using interval arithmetics.
\begin{algorithm}
\hrule height 1pt \vspace{-0.2em}
	\caption{\label{alg::univariate}\small Composition rule of interval superposition arithmetic}
\vspace{0.5em}
\hrule\vspace{0.2em}
	\footnotesize
\smallskip
			\textbf{Input:} Matrix $A$ parameterizing $F_{h,X}$ and an atom
operation $\alpha$.\\[0.1cm]
			\textbf{Main Steps:}\\[-0.3cm]
				\begin{enumerate}[wide]
					\item Choose, for all $i \in \{ 1, \ldots, n_x \}$, central points $a_i \in \mathbb R$
satisfying
					$$L(A_i)  \, \leq \, a_i \, \leq \, U(A_i)  \quad \text{and set} \quad \omega = \sum_{i=1}^{n_x} a_i \; .$$

					\item Choose a suitable remainder bound $r_\alpha(A) \geq 0$ such that
					\begin{eqnarray}
					\left| \sum_{i=1}^{n_x} \alpha( \omega + \delta_i ) - (n_x-1)\alpha(\omega) - \alpha \left( \omega + \sum_{i=1}^{n_x} \delta_i \right) \right| \; \leq \; r_\alpha(A) \notag
					\end{eqnarray}
					for all $\delta \in \mathbb R^{n_x}$ with $\forall i \in \{ 1, \ldots, n_x \}, \; \; L(A_i) \leq a_i + \delta_i \leq U(A_i)$.

					\item Compute the interval valued coefficients
					\[
					C_i^j = \alpha \left(\omega - a_i + A_i^j \right) - \frac{n_x-1}{n_x} \alpha(\omega) \; .
					\]
					for all $i \in \{ 1, \ldots, n_x\}$ and all $j \in \{ 1, \ldots, N \}$, where $\alpha \left(\omega - a_i + A_i^j \right)$ is evaluated in interval arithmetic.
				\item Set $C_k^j \leftarrow C_k^j + r_{\alpha}(A) \cdot [-1,1] $ for all $j \in \{ 1, \ldots, N \}$ with
       \[
        k \in \operatorname*{argmax}_{i\in\{1,\ldots,n_x\}} \sum^{N}_{j=1} \overline{A}_{i}^{j} - \underline{A}_{i}^{j}
       \].
				\end{enumerate}
			\textbf{Output:} Matrix $C\in\mathbb{I}^{n_x\times N}$ parameterizing
      $\Gamma(\cdot,C,X)$ for $\alpha \circ g$.
\vspace{0.5em}
\hrule height 1pt
\end{algorithm}

\begin{algorithm}
\hrule height 1pt \vspace{-0.2em}
	\caption{\label{alg::product}\small Product rule of interval superposition arithmetic}
\vspace{0.5em}
\hrule\vspace{0.2em}
	\footnotesize
	\smallskip
		\textbf{Input:} Matrices $A$ and $B$ parameterizing $F_{h,X}$ and $F_{g,X}$.\\[0.1cm]
		\textbf{Main Steps:}\\[-0.3cm]
			\begin{enumerate}[wide]
				\item Compute the central points, $\forall i \in \{ 1, \ldots, n_x \}$
        \[
        a_i = \frac{\mathrm{U}\left( A_i \right) + \mathrm{L}\left( A_i \right)}{2}\quad
        \text{and} \quad
        b_i = \frac{\mathrm{U}\left( B_i \right) + \mathrm{L}\left( B_i \right)}{2}\,
        \]
				then set
        \[
        a = \sum_{i=1}^{n_x} a_i \; , \; \; b = \sum_{i=1}^n b_i \; , \;\;
        c = \sum_{i=1}^{n_x} a_i b_i \;, \;\;\text{and} \;\;
        \omega = \frac{ab - c}{n_x}\; .
        \]
				\item Compute
				$\rho_i(A) = \frac{U(A_i)-L(A_i)}{2}$ and $\rho_i(B) = \frac{U(B_i)-L(B_i)}{2}$
				for all $i \in \{ 1, \ldots, n_x \}$ as well as the associated remainder bound
				\[
				R(A,B) = \left( \sum_{i=1}^{n_x} \rho_i(A) \right) \left( \sum_{i=1}^{n_x} \rho_i(B) \right) - \sum_{i=1}^{n_x} \rho_i(A)\rho_i(B) \; .
				\]

				\item Compute, for each $i \in \{ 1, \ldots, n_x\}$ and all $j \in \{ 1, \ldots, N \}$
				\[
				C_{i}^j = \left( A_{i}^j + a - a_i \right) \left( B_{i}^j + b - b_i \right) - \left( a - a_i \right)\left(b - b_i \right) - \omega\;.
				\]

				\item Set $C_k^j \leftarrow C_k^j + R(A,B) \cdot [-1,1] $ for all $j \in \{ 1, \ldots, N \}$ with
       \[
        k \in \operatorname*{argmax}_{i\in\{1,\ldots,n_x\}} \sum^{N}_{j=1} \overline{A}_{i}^{j} - \underline{A}_{i}^{j}
       \].
				\end{enumerate}
			\textbf{Output:} Matrix $C\in\mathbb{I}^{n_x\times N}$ parameterizing
      $\Gamma(\cdot,C,X)$, for $g*h$.
\vspace{0.5em}
\hrule height 1pt
\end{algorithm}

\begin{theorem}\label{thm::univariate}
Let $\Gamma(x,A,X)$ and be an ISM of $g$ on $X$. If the matrix
$C\in\mathbb{I}^{n_x\times N}$ is computed using Algorithm~\ref{alg::univariate},
then $\Gamma(x,C,X)$ is an ISM of $\alpha \circ g$ on $X$.
\end{theorem}
\begin{pf}
For the first statement, let $x\in X$ be an arbitrary point. Since $\Gamma(x,A,X) =
\sum^{n_x}_{i=1}\sum^{N}_{j=1}A^{j}_i\varphi^{j}_{i}(x)$ is an ISM of $g$,
there exists a sequence $j_1,\ldots j_{n_x}\in\{1,\dots,N\}$ and points
$y_{i}\in A^{j_i}_{i}$ satisfying $g(x)=\sum^{n_x}_{i=1}y_i$. Let
$\delta_i = y_i -a_i$, with $\omega$ defined as in Algorithm~\ref{alg::univariate}
one can write
\begin{equation*}
\small
\begin{aligned}
\alpha(g(x)) &= \hphantom{{}+{}} \alpha\left( \omega+\sum^{n_x}_{i=1}\delta_i \right)\\
&=\hphantom{{}+{}} \sum^{n_x}_{i=1}\left(\alpha(\omega+\delta_i)
-\frac{n_x-1}{n_x}\alpha(\omega)\right)\\
&\hphantom{{}={}}-\underbrace{\left( \sum^{n_x}_{i=1}\alpha(\omega+\delta_i)-(n_x-1)\alpha)\omega
-\alpha\left(\omega+\sum^{n_x}_{i=1}\delta_i\right) \right)}_{r_{\alpha}(A)[-1,1]}\,.
\end{aligned}
\end{equation*}
Since $\delta_i\in A_{i}^{j_i}-a_i$, we have $\alpha\left(\omega-a_i+A^{j_i}_i\right)$
and
\begin{equation*}
\small
\begin{aligned}
\alpha(g(x)) &\in \sum^{n_x}_{i=1}\left(\alpha\left(\omega-a_i+A^{j_i}_i\right)
-\frac{n_x-1}{n_x}\alpha(\omega)\right) + r_{\alpha}(A)[-1,1]\\
&=\sum^{n_x}_{i=1} C_{i}^{j_i}\;,
\end{aligned}
\end{equation*}
which implies the statement of the theorem.\hfill\hfill\qed
\end{pf}

\begin{theorem}\label{thm::product}
Let $\Gamma(x,A,X)$ and $\Gamma(x,B,X)$ be ISMs of $g$ and $h$, respectively,
on $X$. If $C\in\mathbb{I}^{n_x\times N}$ is computed using
Algorithm~\ref{alg::product}, then  $\Gamma(x,C,X)$ is an ISM of
$g * h$ on $X$.
\end{theorem}

A proof of Thm.~\ref{thm::product} proceeds along the same lines as the proof of
Thm.~\ref{thm::univariate} and its omitted for the sake of brevity.

The construction of remainder bounds and central points used in
Algorithm~\ref{alg::univariate} exploits globally valid algebraic properties,
called addition theorems, of common univariate operations. As an example, for the
exponential function, the addition theorems
$e^{\omega+\delta_i} = e^{\omega}e^{\delta_i}$ and
$e^{\omega+\sum^{n_x}_{i=1}\delta_i} = e^{\omega}\prod^{n_x}_{i=1}e^{\delta_i}$,
hold globally over the real numbers. Letting
$t_i = e^{\delta_i}-1$, $r_{\alpha}(A)$ can be constructed  by bounding
the left-hand side of the expression in Step 2) of Algorithm~\ref{alg::univariate}.
This yields the expression
\begin{equation*}
e^{\omega}\left| \sum^{n_x}_{i=1}t_i+1-\prod^{n_x}_{i=1}(1+t_i)\right|
\leq e^{\omega}\left(  \prod^{n_x}_{i=1}(1+s_i) - \sum^{n_x}_{i=1} s_i-1 \right)
\end{equation*}
with $s_{i}=\max\left\{ e^{U(A_i)-a_i}-1,1-e^{L(A_i)-a_1} \right\}$. Choosing
$a_i=\log\left(\frac{1}{2}\left( e^{U(A_{i})} + e^{L(A_{i})}\right)\right)$,
minimizes
\begin{equation*}
s_{i} = \frac{e^{U(A_i)}-e^{L(A_i)}}{e^{U(A_i)}+e^{L(A_i)}}\;.
\end{equation*}
The technical derivations for the remainder bounds $r_{\alpha}(A)$ and the central
points $a_i$ for other atom operations can be found in~\citep{Zha2016}.

The final ingredient for an arithmetic of interval superposition models is the
construction of a (trivial) ISM for the input variables $x_i$. As each variable
is independent of the rest, the coefficients can be set as $A^{j}_{k} = 0$ for
all $k\neq i$ and all $j\in\{1,\ldots,N\}$. The $i$th row  of $A$ is then
initialized as $A^{j}_{i}=X^{j}_{i}$ for each $j\in\{1,\ldots,N\}$.

Notice that composition rules in ISA have a computational
complexity of ${\rm \bf O}(n_xN)$.
Unfortunately, global statements regarding the model accuracy can only be
given whenever $f$ is separable, i.e.
$f(x) = \sum^{n_x}f_{j}(x_j)$ for some factorable functions $f_{1},\ldots,f_{n_x}$.
In this case, the error is of order ${\rm \bf O}\left(\frac{1}{N}\right)$ over all
bounded domains $X\subset\mathbb{I}^{n_x}$.


\section{ISA-Based Set-inversion Algorithm}
\label{sec::SIVISA}

This section proposes a novel search strategy based on ISA for addresing GPE. It
has as its core computing the intersection of an ISM with an interval.

Consider an ISM, of the function $f$ over $X$, parameterized by
$A\in\mathbb{I}^{n_x\times N}$. The direct way of computing the intersection
between this ISM and $Y=[\underline{y},\overline{y}]$ is to compute the value of
the ISM at each interval $X_{j_1} \times \ldots \times X_{j_N}$ in the partition
of $X$. This requires computing all possible superpositions of coefficients
$A^{i}_{j}$. Such approach, while straightforward, is unfortunately not efficient
since its computational complexity is ${\rm \bf O}\left( N^{n_x} \right)$.

As it turns out, computing an over approximation of the desired intersection can
be done by testing only certain selected combinations. The proposed approach,
requires sorting the components $A^{j}_{i}=[\underline{A}^{j}_{i},
\overline{A}^{j}_{i}]$ of the rows $A_i$ of the matrix $A$ in both decreasing
and increasing orders. The corresponding permutations are denoted by the functions
$\overline{\pi}_{i},\underline{\pi}_{i}:\{1,\ldots,N\}\to\{1,\ldots,N\}$ satisfying
$$\overline{A}^{\overline{\pi}_{i}(1)}_i \geq \overline{A}^{\overline{\pi}_{i}(2)}_i \geq \ldots \geq \overline{A}^{\overline{\pi}_{i}(N)}_i$$
and
$$\underline{A}^{\underline{\pi}_{i}(1)}_i \leq \underline{A}^{\underline{\pi}_{i}(2)}_i \leq \ldots \leq \underline{A}^{\underline{\pi}_{i}(N)}_i \; .$$
In the following, we use the shorthand $\overline \Pi = ( \overline \pi_1, \ldots, \overline \pi_{n_x} )$ and $\underline \Pi = ( \underline \pi_1, \ldots, \underline \pi_{n_x} )$. The main pre-processing step for computing a set
inversion is outlined in Algorithm~\ref{alg:Intersect}.

\begin{algorithm}
\hrule height 1pt \vspace{-0.2em}
	\caption{\small Intersection of a superposition model with an interval}
\vspace{0.5em}
\hrule\vspace{0.2em}
	\label{alg:Intersect}
	\footnotesize
	\bigskip
		\textbf{Input:} Parameters $A$ and $X$ of the input model and an interval $Y$ \\[0.2cm]
		\textbf{Main Step:}\\[-0.3cm]
			\begin{enumerate}[wide]
				\setlength{\itemsep}{2pt}
				\item Sort each $A_i$ to obtain the permutations $\underline \Pi$ and $\overline \Pi$.
				\item Choose a finite number $n_J$ of intervals $\underline J_k = [0, \underline j_k]$ with index vectors $\underline j_k \in \{1, \ldots, N \}^{n_x}$ such that
				$$\forall k\in \{1,\dots,n_J\}, \quad \sum_{i=1}^{n_x} \underline A_i^{{\underline \pi}_i((\underline j_k)_i)} \leq \underline y$$
				\item Choose a finite number $n_J$ of intervals $\overline J_k = [0, \overline j_k]$ with index vectors $\overline j_k \in \{1, \ldots, N \}^{n_x}$ such that
				$$\forall k\in \{1,\dots,n_J\}, \quad \sum_{i=1}^{n_x} \overline A_i^{{\overline \pi}_i((\overline j_k)_i)} \geq \overline y$$
			\end{enumerate}
		\textbf{Output:} Permutations $\underline \Pi,\overline \Pi$ and intervals
$\underline J = (\underline J_1, \ldots, \underline J_{n_J})$,
$\overline J = (\overline J_1, \ldots, \overline J_{n_J})$. \\
\vspace{0.5em}
\hrule height 1pt
\end{algorithm}

\begin{theorem}
\label{thm::Intersect}
Let $\underline \Pi,\overline \Pi$ and $\underline J = (\underline J_1, \ldots, \underline J_{n_J})$, $\overline J = (\overline J_1, \ldots, \overline J_{n_J})$ be computed by Algorithm~\ref{alg:Intersect}. Define
\begin{equation}
 \underline \Xi = \bigcup_{k \in \{ 1, \ldots, n_J\}} \bigcup_{j \in \underline J_k}  \Xi_1^{{\underline \pi}_1(j_1)} \times \ldots \times \Xi_{n_x}^{{\underline \pi}_{n_x}(j_{n_x})}
\end{equation}
and
\begin{equation}
 \overline \Xi = \bigcup_{k \in \{ 1, \ldots, n_J\}} \bigcup_{j \in \overline J_k}  \Xi_1^{{\overline \pi}_1(j_1)} \times \ldots \times \Xi_{n_x}^{{\overline \pi}_{n_x}(j_{n_x})}
\end{equation}
with $\Xi_i^j = [ \underline x_i + (j-1) h_i, \underline x_i + j h_i]$ and $h_i = \frac{\overline x_i - \underline x_i}{N}$. Then,
\begin{equation}
 X \setminus \left( \underline \Xi \cup \overline \Xi \right) \supseteq \mathbb X_\mathrm{int} \cup \mathbb X_\mathrm{bnd} \; .
\end{equation}
\end{theorem}

\begin{pf}
By construction, the function $f$ takes values larger than $\overline y$ on all interval boxes $\Xi_1^{{\underline \pi}_1(j_1)} \times \ldots \times \Xi_{n_x}^{{\underline \pi}_{n_x}(j_{n_x})}$ for any $j \in \overline J_{k}$. Similarly, $f$ takes smaller values than $\underline y$ on all intervals $\Xi_1^{{\overline \pi}_1(j_1)} \times \ldots \times \Xi_{n_x}^{{\overline \pi}_{n_x}(j_{n_x})}$ for any $j \in \underline J_{k}$. Consequently, the union of all of these boxes cannot possibly contain a point of $\mathbb X_\mathrm{int} \cup \mathbb X_\mathrm{bnd}$, which is the statement of the theorem.\hfill\hfill\qed
\end{pf}

Theorem~\ref{thm::Intersect} provides a constructive procedure for finding the desired outer approximation of the set $\mathbb X_\mathrm{int} \cup \mathbb X_\mathrm{bnd}$. Notice that the computational complexity of Algorithm~\ref{alg:Intersect} is of order $O( n_x N \log(N) )$, because we need to sort the intervals along all coordinate directions. The associated storage complexity is of order $O( n_x N )$. Finally, we have to keep in mind, however, that computing and storing the sets $\underline \Xi$ and $\overline \Xi$ is expensive in general, as these sets may be composed of an exponentially large amount of sub-intervals. Nevertheless, it is not necessary to store these sets explicitly as long as we store the permutation matrices $\underline \Pi$ and $\overline \Pi$ as well as the boxes $\underline J$ and $\overline J$, which uniquely represent the set $X \setminus \left( \underline \Xi \cup \overline \Xi \right)$.

Notice that there are various heuristics possible for refining the above procedure.
However, the corresponding methods are analogous to the implementation in SIVIA
and based on state-of-the-art branching techniques. Thus, the proposed technique
based on Algorithm 3 can be embedded in an exhaustive search procedure, if one
wishes to approximate the set $\mathbb X_\mathrm{int} \cup \mathbb X_\mathrm{bnd}$
with any given accuracy.


\section{Numerical Examples}
\label{sec::case}

This section illustrates some of the benefits of ISA as a bounding method for
the range of factorable functions, as well as its application to GPE.
Algorithms~1,~2, and a set-inversion algorithm based on Algorithm~3 were
implemented in the programming language \texttt{Julia}. For comparison, a basic
SIVIA algorithm was also implemented in \texttt{Julia}. The termination for both
algorithms was based on~\eqref{eq::maxdiam}. All results
were obtained on an Intel Xeon CPU X5660 with 2.80GHz and 16GB RAM.

\subsection{Bounding a nonlinear function: ISMs vs TMs}

Consider the nonlinear factorable function
\begin{equation*}
 f(x) = e^{\sin(x_1)+\sin(x_2)\cos(x_2)}
\end{equation*}
over the domain $X = [0,1]\times[0,\overline{x}_{2}]$\;. Here,
$\overline{x}_{2}\in[0.1,20]$ denotes a parameter which controlling the diameter
of the domain. In order to measure the quality of an arithmetic, we used the
Hausdorff distance between the range of $f$, $f(X)={f(x)|x\in X}$, and an
enclosure set $\overline Y \supseteq f(X)$. This distance is given by
\begin{equation*}
d_{\rm H}(f(X),\overline{Y}) = \max_{y\in \overline Y} \min_{x\in f(X)} || x-y||_{\infty}\;.
\end{equation*}

\begin{figure}
\centering
\begin{overpic}[width=0.7\columnwidth]{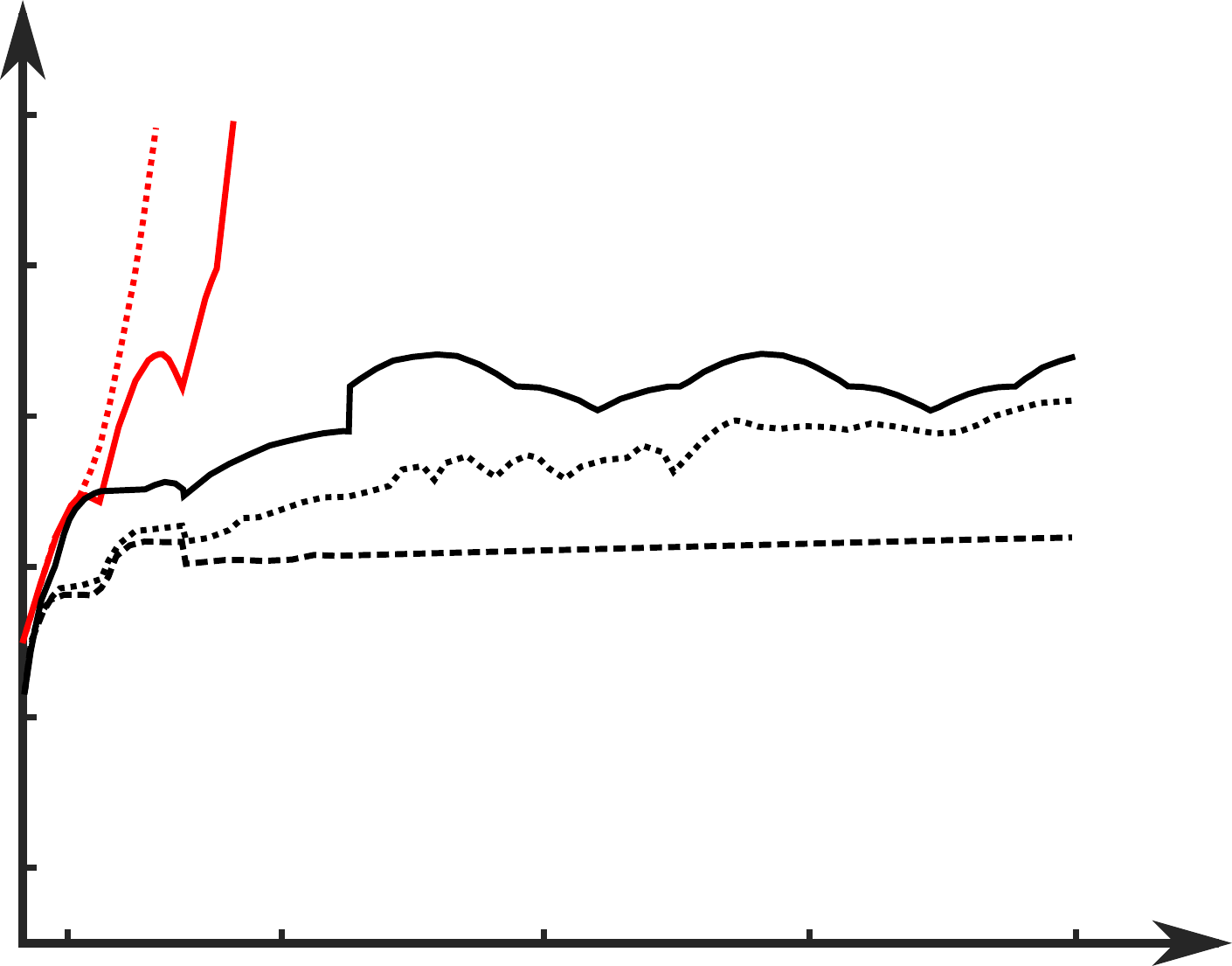}
\put(5,74){$d_{\rm H}\left(f(X),\overline Y \right)$}
\put(94,6){$\overline x_2$}

\put(88,49){$N = 1$}
\put(88,43){$N = 10$}
\put(88,34){$N = 100$}

\put(20,59){\rotatebox{88}{\color{red}TM1}}
\put(6,59){\rotatebox{84}{\color{red}TM2}}

\put(4,-4){$1$}
\put(22,-4){$5$}
\put(42,-4){$10$}
\put(64,-4){$15$}
\put(85,-4){$20$}

\put(-10,6){$10^{-2}$}
\put(-9.9,18){$10^{-1}$}
\put(-2.5,30){$1$}
\put(-7,42.5){$10^1$}
\put(-7.2,55){$10^2$}
\put(-7.2,67){$10^3$}
\end{overpic}
\vspace{0.5em}
\caption{\small \label{fig::enclosure}Overestimation of enclosure
sets with respect to the domain size. The plot compares enclosures based on TMs of
orders~1 (solid red) and~2 (dotted red) as well as ISMs with $N=1$ (solid black),
$N=10$ (dotted black), and $N=100$ (dashed black).}
\end{figure}
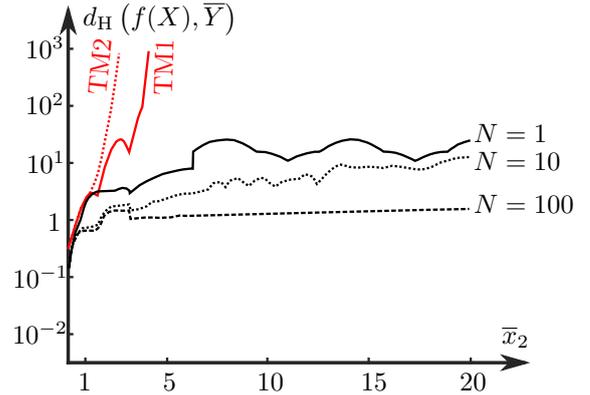

Figure~\ref{fig::enclosure} shows the overestimation of enclosures
in the form of Taylor models of orders 1 and 2 as well as interval superposition
models with $N=1$, $N=10$, and $N=100$ as a function of the domain parameter
$\overline{x}_2$. Although the Hausdorff distance between $f(X)$ and $\overline{Y}$
does not increase monotonically with $\overline{x}_2$, the rough
trend observed on the plot is that the overestimation increases with the
size of the domain. Furthermore, the plot shows that interval superposition models
outperform Taylor models over large domains. One aspect that is not shown in the
figure is that over small domains, e.g. over $[0,10^{-1}]^2$, enclosures based on
Taylor models outperform those constructed using interval superposition
arithmetics.

\subsection{Guaranteed parameter estimation via ISMs}
{\color{black}

We consider a reaction system given :
\begin{equation}
\label{eq::odes}
\begin{alignedat}{2}
\dot{z_1}(t) &= - (x_1 + x_3) z_1(t) + x_2 z_2(t), \quad &&z_1(0) = 1, \\
\dot{z_2}(t) &= x_1 z_1(t) - x_2 z_2(t), &&z_2(0) = 0\;,
\end{alignedat}
\end{equation}
with $y(t) = z_2(2)$~\citep{pau15-ima}. The output variable, can be
represented as the factorable function
\begin{equation*}
y(t) = \frac{e^{\frac{-t\rho}{2}}x_1 (e^{\frac{t\sigma}{2}} - e^{\frac{ - t\sigma}{2}})}{\sigma}
\end{equation*}
with $\sigma = \sqrt{x_1^2 + x_2^2 + x_3^2 + 2x_1p_2 + 2x_1x_3 - 2x_2x_3}$ as well
as $\rho = x_1 + x_2 + x_3$. In the following, we fix $x_3 = 0.35$ and consider
$n_m = 15$ measurements corresponding to the time instants $t_i=1,2,\ldots,15$.
Process measurements were obtained by simulating~\eqref{eq::odes} with
$x = (0.6,0.15,0.35)^T$, rounding to the second significant digit. Measurement
errors of $\pm 10^{-3}$ were added to these values.

The performance of the proposed GPE algorithm using ISA was tested against a
standard SIVIA. We have interval superposition models with $N=2,\,10,\,20$.
Figure~\ref{fig::results} shows a summary of the results of the GPE algorithm
using ISMs with $N=2$. The left plot, shows an approximation of the set $X_e$.
The plot shows the inner partition (in red) for $\epsilon=10^{-5}$ and the
boundary partitions for $\epsilon=10^{-4}$ (light blue) and $\epsilon-10^{-5}$
(dark blue). The central and right plots show, respectively, a comparison of the
number of iterations and CPU time against the tolerance $\epsilon$---for SIVIA
(solid red line) and ISM-based set-inversion with $N=2$ (solid black line),
$N=10$ (dotted black line), and $N=20$ (dashed black line). In terms of the
number of iterations and number of boundary boxes (not shown), ISM-based
set-inversion (for all $N$) outperforms SIVIA. This is due to the fact that
ISA is able to detect and exploit structures in the factorable function to
remove redundant boxes. On the contrary, with respect to the CPU time, SIVIA
outperforms the proposed algorithm. This can be traced back to the fact that
the the cost per iteration is larger for ISA. Furthermore, the implementation
is still at prototype stage and requires further refinement in terms of computing
remainder bounds and memory management in the algorithms.

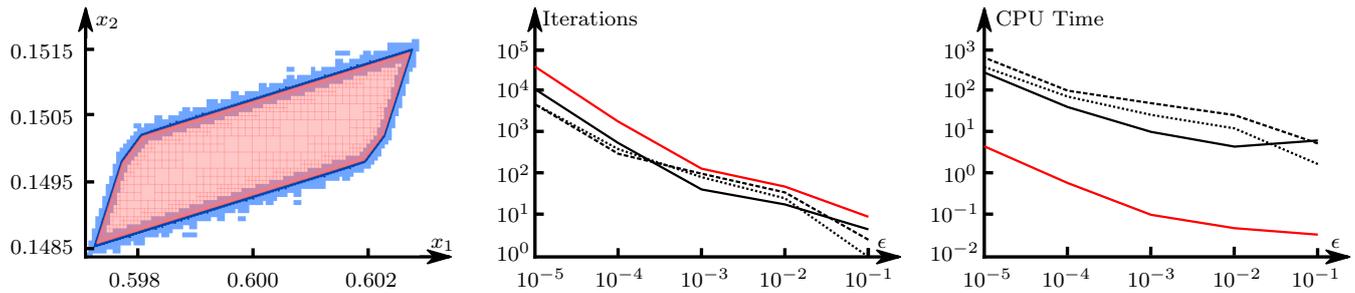
\begin{figure*}
\centering
\begin{overpic}[width=0.95\textwidth]{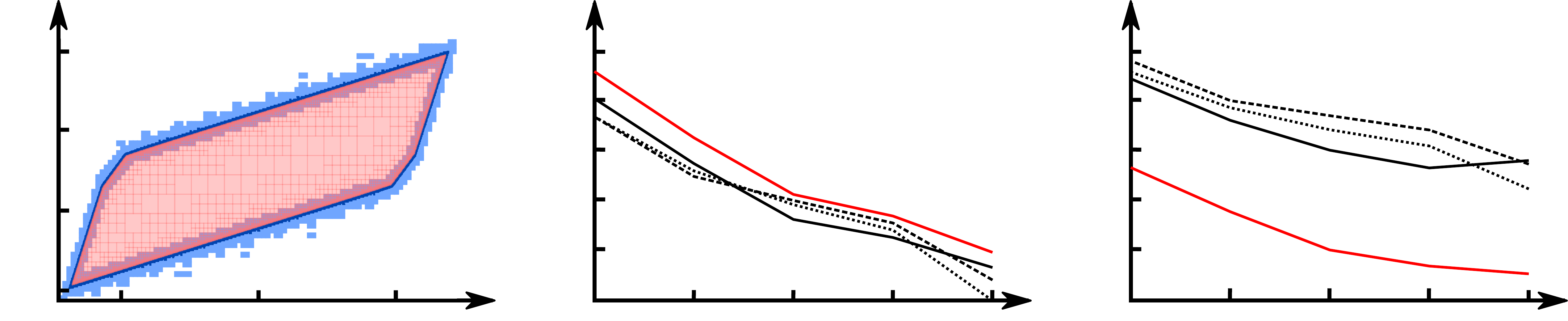}

\put(30,2){\small$x_1$}
\put(5.5,-1){\small$0.598$}
\put(14.5,-1){\small$0.600$}
\put(23.5,-1){\small$0.602$}

\put(4.5,19){\small$x_2$}
\put(-2,1.5){\small$0.1485$}
\put(-2,6.5){\small$0.1495$}
\put(-2,11.5){\small$0.1505$}
\put(-2,16.5){\small$0.1515$}

\put(64,2){\small$\epsilon$}
\put(36.5,-1){\small$10^{-5}$}
\put(42.5,-1){\small$10^{-4}$}
\put(49,-1){\small$10^{-3}$}
\put(55,-1){\small$10^{-2}$}
\put(61.5,-1){\small$10^{-1}$}

\put(38.5,19){\small Iterations}
\put(35,1){\small$10^{0}$}
\put(35,4){\small$10^{1}$}
\put(35,7){\small$10^{2}$}
\put(35,10.5){\small$10^{3}$}
\put(35,13.5){\small$10^{4}$}
\put(35,16.5){\small$10^{5}$}

\put(98.5,2){\small$\epsilon$}
\put(70.5,-1){\small$10^{-5}$}
\put(76.5,-1){\small$10^{-4}$}
\put(83.0,-1){\small$10^{-3}$}
\put(89.5,-1){\small$10^{-2}$}
\put(95.5,-1){\small$10^{-1}$}

\put(73,19){\small CPU Time}
\put(68,1){\small$10^{-2}$}
\put(68,4){\small$10^{-1}$}
\put(69,7){\small$10^{0}$}
\put(69,10.5){\small$10^{1}$}
\put(69,13.5){\small$10^{2}$}
\put(69,16.5){\small$10^{3}$}

\end{overpic}
\vspace{0.5em}
\caption{\small \label{fig::results} Results for the GPE problem. Left: Parameter
inner partition (in red) for $\epsilon=10^{-5}$ and the boundary partitions for
$\epsilon=10^{-4}$ (light blue) and $\epsilon-10^{-5}$ (dark blue). Center:
Number of iterations vs. diameter of boundary partition. Right: CPU time vs. diameter
of boundary partition. Center and right plots show results for SIVIA
(solid red line) and ISM-based set-inversion with $N=2$ (solid black line),
$N=10$ (dotted black line), and $N=20$ (dashed black line).}
\end{figure*}

}

\section{Conclusion}
\label{sec::conclusion}
This paper presented Interval superposition arithmetics, a novel set-arithmetic
for computing enclosures of the image set of factorable functions and its use in
guaranteed parameter estimation. The main advantage of ISA is its polynomial
storage and computational complexity. {\color{black}The core routine behind the
proposed GPE method is the intersection of an interval superposition model and
an interval. Although the proposed intersection routine has a computational
complexity of order $O(n_x N \log(N))$, computing an arbitrarily accurate
approximation of the parameter set requires exponential run time. Our numerical
examples illustrate the advantages of ISA over other set arithmetics when
constructing enclosures for factorable functions---particularly over large domains.
We have also shown how the proposed technique can be used to solve a GPE problem.
Although the number of iterations is reduced when using ISA, the overal CPU time
is larger than SIVIA. This suggest that, although ISA can improve certain aspects
of GPE algorithms, there is still much room for improvement. Improved ISA-based
algorithms for constructing approximations of inverse-image sets in polynomial
run-time will be investigated in future work.}

\section*{Acknowledgments}{\small
This work was supported by: National Science Foundation China (NSFC), Grant~61473185;
ShanghaiTech University, Grant~F-0203-14-012; and Slovak Research and
Development Agency, project APVV SK-CN-2015-0016 ``CN-SK cooperation:
Verified Estimation and Control of Chemical Processes''. RP also acknowledges:
Slovak Research and Development Agency, project APVV 15-0007; and European
Commission, grant agreement 790017 (GuEst).
}
\bibliography{references}

\end{document}